\documentclass[11pt, a4paper]{amsart}

\usepackage[abbrev]{amsrefs}

\usepackage{cmdtrack}

\usepackage[latin1]{inputenc}
\usepackage{enumerate}

\newcommand{\R}{\mathbb R}
\newcommand{\N}{\mathbb N}
\newcommand{\cM}{\mathcal{M}}
\newcommand{\e}{\mathrm{e}}
\newcommand{\dif}{\, \mathrm{d}}
\newcommand{\diam}{\operatorname{diam}}
\newcommand{\sign}{\operatorname{sign}}
\DeclareMathOperator{\capt}{cap}
\newtheorem{teorema}{Theorem}[section]
\newtheorem{proposicao}[teorema]{Proposition}
\newtheorem{corolario}[teorema]{Corollary}

\title{Limit solutions of the Chern-Simons equation}

\author{Augusto C. Ponce}
\address{
Augusto C. Ponce\hfill\break\indent
 Universit{\'e} catholique de Louvain\hfill\break\indent
 Institut de Recherche en Math{\'e}matique et Physique\hfill\break\indent
 Chemin du cyclotron 2, bte L7.01.02\hfill\break\indent
1348 Louvain-la-Neuve\hfill\break\indent
Belgium}
\email{Augusto.Ponce@uclouvain.be}

\author{Adilson E. Presoto}
\address{
Adilson E. Presoto\hfill\break\indent
Universidade Estadual de Campinas\hfill\break\indent
IMECC, Departamento de Matem\'atica\hfill\break\indent
Rua S\'ergio Buarque de Holanda, 651\hfill\break\indent
Campinas, SP, 13083-859\hfill\break\indent
Brazil}
\email{ra057746@ime.unicamp.br}

\begin{document}

\begin{abstract}
Given a bounded domain $\Omega$ in $\mathbb{R}^2$, we investigate the scalar Chern-Simons equation 
\[
-\Delta u + \e^u(\e^u-1) =\mu \quad \text{in \(\Omega\),}
\]
in cases where there is no solution for a given nonnegative finite measure  $\mu$. 
Approximating $\mu$ by a sequence \((\mu_n)_{n \in \N}\) of nonnegative \(L^1\) functions or finite measures for which this equation has a solution, we show that the sequence of solutions \((u_n)_{n \in \N}\) of the Dirichlet problem converges to the solution with largest possible datum \(\mu^\#  \le \mu\)  and we derive an explicit formula of \(\mu^\#\) in terms of \(\mu\).
The counterpart for the Chern-Simons system with datum $(\mu, \nu)$ behaves differently and the conclusion depends on how much the measures \(\mu\) and \(\nu\) charge singletons.
\end{abstract}

\keywords{elliptic system, exponential nonlinearity, scalar Chern-Simons equation, Chern-Simons system, Radon measure}

\subjclass[2010]{Primary 35R06; Secondary 35J25, 35J57} 

\maketitle

\section{Introduction and main results}

In this paper we investigate a question concerning convergence and stability 
of solutions of the scalar Chern-Simons problem
 \begin{equation}
 \label{eq01}
  \left\{
   \begin{alignedat}{2}
    -\Delta u + \e^u(\e^u - 1) & = \mu && \quad \mbox{in}\; \Omega,\\
                u &=0  && \quad \mbox{on}\; \partial \Omega,
   \end{alignedat}
  \right.
 \end{equation}
where \(\Omega \subset \R^2\) is a smooth bounded domain and \(\mu\) is a finite Borel measure --- equivalently a Radon measure --- in \(\Omega\). 
By a solution of \eqref{eq01}, we mean a function \(u \in W_0^{1, 1}(\Omega)\) such that $\e^u(\e^u-1)\in L^1(\Omega)$ and satisfying the equation in the sense of distributions.

Using for instance a minimization argument in \(W_0^{1, 2}(\Omega)\), one shows that the scalar Chern-Simons equation always has a solution with datum \(\mu \in L^p(\Omega)\) for any \(1 < p \le \infty\) \cite{Ponce}*{Chapter~2}. 
Existence in the case of datum \(\mu \in L^1(\Omega)\) can be obtained by approximation using \(L^\infty\) data~\citelist{\cite{BS}*{Corollary~12} \cite{Ponce}*{Chapter~3}}.

The case of nonlinear Dirichlet problems with measure data is more subtle.
This issue has been discovered by Bénilan and Brezis~\cites{BenBre:04,Bre:80,Bre:82} in a pioneering work concerning polynomial nonlinearities in dimension greater than \(2\).

The case of exponential nonlinearities in dimension \(2\) has been investigated by Vázquez~\cite{V}. 
For instance, if \(\mu = \alpha \delta_a\) for some \(a \in \Omega\), then for every \(\alpha > 2\pi\) the Dirichlet problem~\eqref{eq01} has no solution with datum \(\mu\).
The counterexample above gives the only possible obstruction in the case of exponential nonlinearities: \(\mu\) is a good measure --- that is the Dirichlet problem \eqref{eq01} has a solution --- if and only if for every \(x \in \Omega\), 
\[
\mu(\{x\}) \le 2\pi.
\]

We want to understand what happens when one forces the Dirichlet problem to have a solution when no solution is available.
For instance, if \(\mu\) is a measure for which \eqref{eq01} has no solution, then one could approximate \(\mu\) by a sequence \((\rho_n * \mu)_{n \in \N}\) of convolutions of \(\mu\) --- for which we know the Dirichlet problem has a solution --- and then investigate the limit of the sequence of solutions \((u_n)_{n \in \N}\).

This program has been proposed and implemented by Brezis, Marcus and Ponce~\cite{BPM} in the case where \(\mu\) is approximated via convolution. 
They have proved that for any sequence of nonnegative mollifiers \((\rho_n)_{n \in \N}\), if \(u_n\) satisfies
 \begin{equation*}
  \left\{
   \begin{alignedat}{2}
    -\Delta u_n + \e^{u_n}(\e^{u_n} - 1) & = \rho_n *\mu && \quad \mbox{in } \Omega,\\
     u_n &=0 && \quad \mbox{on }\partial \Omega,
   \end{alignedat}
  \right.
 \end{equation*}
then the sequence \((u_n)_{n \in \N}\) converges in \(L^1(\Omega)\) to the largest subsolution \(u^*\) of the scalar Chern-Simons problem with datum \(\mu\)~\cite{BPM}*{Theorem~4.11}.

The result in \cite{BPM} concerns more general convex nonlinearities and holds in any dimension, but strongly relies on the fact that the approximating sequence \((\rho_n * \mu)_{n \in \N}\) is constructed via convolution of \(\mu\)~\cite{BPM}*{Example~4.1}.

Our first result shows that for the Chern-Simons equation the conclusion of Brezis, Marcus and Ponce is always true regardless of the sequences of functions --- or even measures --- \((\mu_n)_{n \in \N}\) converging to \(\mu\).

\begin{teorema}
\label{teocc}
Let $(\mu_n)_{n \in \N}$ be a nonnegative sequence of measures in \(\Omega\) such that for every \(n \in \N\) and for every \(x \in \Omega\), 
\[
\mu_n(\{x\}) \le 2\pi
\]
and let \(u_n\) satisfy the scalar Chern-Simons problem
\[
  \left\{
   \begin{alignedat}{2}
    -\Delta u_n + \e^{u_n}(\e^{u_n} - 1) & = \mu_n && \quad \mbox{in } \Omega,\\
                u_n &=0  && \quad \mbox{on } \partial \Omega.
   \end{alignedat}
  \right.
\]
If the sequence \((\mu_n)_{n \in \N}\) converges to a measure \(\mu\) in the sense of measures in \(\Omega\), 
then the sequence $(u_n)_{n \in \N}$ converges in \(L^1(\Omega)\) to the solution of the scalar Chern-Simons problem with datum \(\mu^\#\), where \(\mu^\#\) is the largest measure less than or equal to \(\mu\) such that for every \(x \in \Omega\), 
\[
\mu^\#(\{x\}) \le 2\pi. 
\]
\end{teorema}

A sequence \((\mu_n)_{n \in \N}\) converges weakly to \(\mu\) in the sense of measures in \(\Omega\), if for every continuous function \(\zeta : \overline\Omega \to \R\) such that \(\zeta = 0\) on \(\partial\Omega\),
\[
\lim_{n \to \infty}{\int_\Omega \zeta \dif \mu_n} = \int_\Omega \zeta \dif \mu.
\]
We denote this convergence by  $\mu_n\overset{*}\rightharpoonup \mu$ in $\mathcal{M}(\Omega)$, where $\mathcal{M}(\Omega)$ is the vector space of (finite) measures in \(\Omega\) equipped with the norm
\[
\|\mu\|_{\mathcal{M}(\Omega)}
= |\mu|(\Omega)
= \int\limits_\Omega \dif |\mu|.
\]

Applying Theorem~\ref{teocc} we deduce an explicit formula of \(\mu^\#\) in terms of \(\mu\). 
Indeed, if we write \(\mu\) as a sum of nonatomic part \(\overline{\mu}\) and an atomic part
\[
\mu = \overline{\mu} + \sum_{i=0}^\infty {\alpha_i \delta_{a_i}},
\]
where \(\alpha_i \ge 0\) and the points \(a_i\) are distinct, then
\[
\mu^\# = \overline{\mu} + \sum_{i=0}^\infty {\min{\{\alpha_i, 2\pi\}} \delta_{a_i}}.
\]
Since \(\mu\) is a finite measure, there can only be finitely many indices \(i\) such that \(\alpha_i > 2\pi\).
In particular, the measure \(\mu - \mu^\#\) is supported in a finite set and for every \(a \in \Omega\),
\[
\mu^\#(\{a\}) = \min{\{\mu(\{a\}), 2\pi\}}.
\]

We may recover the result of Brezis, Marcus and Ponce using their notion of reduced measure \(\mu^*\).
By definition, the reduced measure is the unique locally finite measure in \(\Omega\) such that
\[
\mu^* = - \Delta u^* + \e^{u^*}(\e^{u^*} - 1)
\]
in the sense of distributions in \(\Omega\), where \(u^*\) is the largest subsolution of the Dirichlet problem~\eqref{eq01}.
The fundamental property of reduced measures~\cite{BPM}*{Theorem~4.1} asserts  that \(\mu^*\) is a (finite) measure in \(\Omega\), \(u^*\) satisfies the Dirichlet problem~\eqref{eq01} with datum \(\mu^*\) and \(\mu^*\) is the largest good measure less than or equal to \(\mu\).
According to Vázquez's result such largest good measure is precisely \(\mu^\#\).
Therefore,
\[
\mu^\# = \mu^*.
\]

\medskip

As an application of the tools we use to prove Theorem~\ref{teocc}, we investigate what happens to the approximation scheme for the Chern-Simons system
\[
\left\{
\begin{alignedat}{2}
    -\Delta u+ \e^{v}(\e^{u}-1)&=\mu && \quad \text{in \(\Omega\)},\\
    -\Delta v+\e^{u}(\e^{v}-1) &=\nu && \quad \text{in \(\Omega\)},\\
    u=v & =0 && \quad \text{on \(\partial \Omega\)}.
   \end{alignedat}
   \right.
 \]

By a result of Lin, Ponce and Yang~\cite{LPY}*{Theorem~1.1}, the system above has a solution for nonnegative measures \(\mu\) and \(\nu\) in \(\Omega\) if and only if for every \(x \in \Omega\),
\[
\mu(\{x\})+\nu(\{x\})\leq 4\pi.
\]

A first result in this direction consists in identifying the sum of the components of the \emph{reduced limit} \((\mu^\#, \nu^\#)\).

\begin{teorema}
\label{teoprinsist}
Let \((\mu_n)_{n \in \N}\) and \((\nu_n)_{n \in \N}\) be sequences of nonnegative measures such that for every \(n \in \N\) and \(x \in \Omega\),
\[
\mu_n(\{x\}) + \nu_n(\{x\}) \le 4\pi
\] 
and let \((u_n, v_n)\) satisfy the Chern-Simons system
\[
 \label{sistcs}
\left\{
 \begin{alignedat}{2}
    -\Delta u_n+\e^{v_n}(\e^{u_n}-1)&=\mu_n && \quad \mbox{in } \Omega,\\
    -\Delta v_n+\e^{u_n}(\e^{v_n}-1)&=\nu_n && \quad \mbox{in } \Omega,\\
    u_n=v_n & =0 && \quad\mbox{on } \partial \Omega.
   \end{alignedat}
   \right.
\]
If the sequences \((\mu_n)_{n \in \N}\) and \((\nu_n)_{n \in \N}\) converge to \(\mu\) and \(\nu\) in the sense of measures in \(\Omega\), and if \((u_n)_{n \in \N}\) and \((v_n)_{n \in \N}\) converge to \(u\) and \(v\) in \(L^1(\Omega)\),
then $(u, v)$ satisfies the Chern-Simons system
with datum \((\mu^{\#}, \nu^\#)\), where 
\begin{enumerate}[\((i)\)]
\item \(0 \le \mu^\# \le \mu\),
\item \(0 \le \nu^\# \le \nu\),
\item \(\mu^\# + \nu^\#\) is the largest measure less than or equal to \(\mu + \nu\) such that for every \(x \in \Omega\),
\[
\mu^\#(\{x\}) + \nu^\#(\{x\}) \le 4\pi.
\]
\end{enumerate}
\end{teorema}

The assumption concerning the convergence of the sequences \((u_n)_{n \in \N}\) and \((v_n)_{n \in \N}\) in $L^1(\Omega)$ is not restrictive since
if the sequences \((\mu_n)_{n \in \N}\) and \((\nu_n)_{n \in \N}\) are bounded \(\cM(\Omega)\), then both \((u_n)_{n \in \N}\) and \((v_n)_{n \in \N}\) are compact in $L^1(\Omega)$.
This is a consequence of a contraction estimate for Dirichlet problems with an absorption nonlinearity and Stampacchia's linear regularity theory~\citelist{\cite{Stampacchia}*{Théorème~9.1} \cite{Ponce}*{Chapter~3 and Chapter~4}}. 

By Theorem~\ref{teoprinsist}, the measure \(\mu^\# + \nu^\#\) can be expressed only in terms of the measure \(\mu + \nu\) as we did in the case of the scalar Chern-Simons problem.
We deduce that the nonatomic parts of \(\mu^\#\) and \(\nu^\#\) coincide with the nonatomic parts of \(\mu\) and \(\nu\),
so we are still left to identify the atomic parts of \(\mu^\#\) and \(\nu^\#\).

This raises the following question:
are the measures \(\mu^\#\) and \(\nu^\#\) independent of the choice of sequences \((\mu_n)_{n \in \N}\) and \((\nu_n)_{n \in \N}\) converging weakly in measure to \(\mu\) and \(\nu\)?
Despite of what happens in the scalar case, in general the answer is negative; see Section~\ref{NRDCSS}.

We have been able to identify two cases where the answer is affirmative.

\begin{teorema}
\label{thmChernSimonsSystemAtomicPart}
Let \(a \in \Omega\). 
Under the assumptions of Theorem~\ref{teoprinsist},
\begin{enumerate}[\((i)\)]
\item if \(\mu(\{a\})= 0\) or \(\nu(\{a\}) = 0\), then
\begin{align*}
  \mu^{\#}
  (\{a\})&=\min{\left\{\mu(\{a\}), 4\pi \right\}}, \\    
  \nu^{\#}(\{a\})&=\min{\left\{\nu(\{a\}), 4\pi \right\}};
\end{align*}
\item if \(\mu(\{a\})\leq 4\pi\) and \(\nu(\{a\})\leq 4\pi\), then
\begin{align*}
  \mu^{\#}
  (\{a\})&=\min{\left\{\mu(\{a\}),\  \mu(\{a\}) - \frac{\mu(\{a\}) + \nu(\{a\}) - 4\pi}{2}\right\}}, \\    
  \nu^{\#}(\{a\})&=\min{\left\{\nu(\{a\}),\ \nu(\{a\}) - \frac{\mu(\{a\}) + \nu(\{a\}) - 4\pi}{2} \right\}}.
\end{align*}
\end{enumerate}
\end{teorema}

Without one of these assumptions on the measures \(\mu\) and \(\nu\),  the reduced limit \((\mu^\#, \nu^\#)\) depends on the choice of sequences \((\mu_n)_{n \in \N}\) and \((\nu_n)_{n \in \N}\).
The proof is based on Theorem~\ref{thmChernSimonsSystemAtomicPart} and on Cantor's diagonal argument.
A sketch of the argument is presented in Section~\ref{NRDCSS} below.

\section{Proof of Theorem~\ref{teocc}}
\label{Thm1}

By a standard property of elliptic equations with absorption term \cite{Ponce}*{Chapter~7}, for every \(n \in \N\),
\begin{equation}
\label{eqEstimationL1}
\|\e^{u_{n}}(\e^{u_{n}}-1)\|_{L^1(\Omega)} 
\le \|\mu_n\|_{\cM(\Omega)}.
\end{equation}
Thus, by the triangle inequality,
\[
\|\Delta u_n\|_{\cM(\Omega)} 
\le 2 \|\mu_n\|_{\cM(\Omega)}.
\]
Since the sequence \((\mu_n)_{n \in \N}\) is bounded in \(\cM(\Omega)\), the sequence \((\Delta u_n)_{n \in \N}\) is also bounded in \(\cM(\Omega)\).
From Stampacchia's linear regularity theory~\citelist{\cite{Stampacchia}*{Théorème~9.1} \cite{Ponce}*{Chapter~3}}, 
the sequence \((u_n)_{n \in \N}\) is bounded in \(W^{1, q}(\Omega)\) for every \(1 \le q < 2\).
By the Rellich-Kondrachov compactness theorem, there exists a subsequence $(u_{n_k})_{k \in \N}$ converging to some function \(u\) in \(L^1(\Omega)\) and a.e.~in \(\Omega\).
By \eqref{eqEstimationL1}, the sequence \((\e^{u_{n}}(\e^{u_{n}} - 1))_{n \in \N}\) is bounded in \(L^1(\Omega)\). 
Passing to a further subsequence if necessary, we may assume that there exists a finite measure $\tau$ in \(\Omega\) such that 
\begin{equation*}
\e^{u_{n_k}}(\e^{u_{n_k}}-1) \stackrel{\ast}{\rightharpoonup} \e^u(\e^u-1) + \tau \quad \text{in \(\cM(\Omega)\).}
\end{equation*}
Thus, $u$ satisfies the scalar Chern-Simons problem
\begin{equation*}
 \left\{
    \begin{alignedat}{2}
     -\Delta u+ \e^u(\e^u-1) &= \mu-\tau  &&  \quad \mbox{in } \Omega,\\
                   u&=0 && \quad \mbox{on } \partial \Omega.
    \end{alignedat}
   \right.
  \end{equation*}

Consider the set  
\[
A = \big\{x\in \Omega : \mu(\{x\})\geq 2\pi \big\}.
\]
Since \(\mu\) is a finite measure, the set \(A\) is finite.
We first prove that $\tau$ is supported in \(A\). 

For this purpose, let \(N(\mu_n)\) be the Newtonian potential generated by \(\mu_n\),
\[
N(\mu_n)(x) 
= \frac{1}{2\pi} \int_\Omega \log{\left( \frac{d}{|x - y|} \right)} \dif \mu_n(y),
\]
where \(d \ge \diam{\Omega}\).
Given \(b \in \Omega\) and \(r > 0\), we first write the Newtonian potential of \(\mu_n\) as
\[
N(\mu_n)
= N(\mu_n\lfloor_{B_r(b)}) + N(\mu_n\lfloor_{\Omega \setminus B_r(b)}).
\] 
Assume for the moment that there exist \(\epsilon > 0\) and \(m \in \N\) such that for every \(n \ge m\),
\begin{equation}\label{estimate}
\mu_n(B_r(b)) \leq 2\pi- \epsilon.
\end{equation}
By the Brezis-Merle inequality~\citelist{\cite{BM}*{Theorem~1} \cite{Ponce}*{Lemma~8.2}}, there exist $p>1$ and $C_1 >0$ such that for every \(n \ge m\),
 \begin{equation*}
  \big\|\e^{2 N(\mu_n\lfloor_{B_r(b)})}\big\|_{L^p(\Omega)} \leq C_1.
 \end{equation*}
Since the functions \(N(\mu_n\lfloor_{\Omega \setminus B_r(b)})\) are harmonic in \(B_r(b)\) and have a uniformly bounded \(L^1\) norm in \(B_r(b)\), the sequence \(\big(N(\mu_n\lfloor_{\Omega \setminus B_r(b)})\big)_{n \in \N}\) is uniformly bounded in \(B_{r/2}(b)\).
We conclude that there exists \(C_2 > 0\) such that for every \(n \ge m\),
\begin{equation}\label{estimate-bis}
  \big\|\e^{2 N(\mu_n)}\big\|_{L^p(B_{r/2}(b))} \leq C_2.
\end{equation}

Note that if \(b \in \Omega \setminus A\), then there exist \(\epsilon > 0\)
and \(r > 0\) satisfying \eqref{estimate}.
Indeed, let \(\overline{\epsilon} > 0\) and \(R > 0\) such that 
\[
\mu(B_R(b)) \le 2\pi - \overline{\epsilon}.
\]
Then, by weak convergence of the sequence \((\mu_n)_{n \in \N}\) \cite{Evans_Gariepy}*{Section~1.9}, property \eqref{estimate} holds for every \(0 < r < R\) and for every \(0 < \epsilon < \overline{\epsilon}\).

Let \(U_n\) be the solution of the linear Dirichlet problem
 \begin{equation}
 \label{eqLinearDP}
  \left\{
   \begin{alignedat}{2}
    -\Delta U_n&=\mu_n && \quad \text{in }\Omega,\\
            U_n&=0     && \quad \text{on }\partial \Omega.
   \end{alignedat}
  \right.
 \end{equation}
By the comparison estimate between the solution \(U_n\) of the linear Dirichlet problem \eqref{eqLinearDP} and the solution \(u_n\) of the nonlinear Dirichlet problem \cite{Ponce}*{Chapter~7}, for every $n\in \mathbb{N}$ we have
\[
u_n\leq U_n \quad  \text{in $\Omega$. }
\]
By the weak maximum principle \cite{Ponce}*{Chapter~5}, \(U_n \le N(\mu_n)\) in \(\Omega\). 
Hence,
\[
u_n\leq N(\mu_n) \quad  \text{in $\Omega$. }
\]
It follows from \eqref{estimate-bis} that the sequence $(\e^{u_{n}}(\e^{u_{n}} - 1))_{n \in \N}$ is uniformly bounded in $L^{p}(B_{r/2}(b))$.
Since $u_{n_k}\rightarrow u$ a.e. in $B_{r/2}(b)$, by Egorov's theorem we get 
\[
\e^{u_{n_k}}(\e^{u_{n_k}} - 1) \rightarrow \e^{u}(\e^{u} - 1) \quad \text{in \(L^{1}(B_{r/2}(b))\).}
\]
We deduce that $\tau=0$ in $B_{r/2}(b)$. 
Since \(b \in \Omega \setminus A\) is arbitrary, we conclude that \(\tau\) is supported in \(A\).

\medskip
If the set \(A\) is empty, the conclusion of the theorem follows with \(\mu^\# = \mu\).
We may assume that \(A\) is nonempty.
Recalling that \(A\) is a finite set, we may write
\[
A = \{x_1, \dots, x_l\},
\]
where the points \(x_i \in \Omega\) are distinct.

Given $i \in \{1, \dots, l\}$, let \(r > 0\) be such that \(B_r(x_i) \cap A = \{x_i\}\). 
For every 
\[
0 \le \alpha < \frac{2 \pi}{\mu(\{x_i\})},
\]
let \(v_k\) be a function satisfying the scalar Chern-Simons problem
\[
  \left\{
   \begin{alignedat}{2}
    -\Delta v_{{k}}+ \e^{v_k}(\e^{v_{{k}}}-1) &= \alpha \mu_{n_k} && \quad \mbox{in } B_r(x_i),\\
     v_{k}&=0 && \quad \mbox{on } \partial B_r(x_i).
   \end{alignedat}
  \right.
\]
The existence of $v_{k}$ follows from \cite{V}*{Theorem~2}; alternatively, one may apply the method of sub and supersolution~\citelist{\cite{MP}*{Corollary~5.4} \cite{Ponce}*{Chapter~6}} with subsolution \(0\) and supersolution \(u_{n_k}\). 
In particular,
\[
0 \le v_k \le u_{n_k} \quad \text{in } \Omega.
\]

Since for every \(x \in B_r(x_i)\), 
\[
\alpha \mu(\{x\}) < 2\pi,
\]
the sequence \((v_k)_{k \in \N}\) converges in \(L^1(\Omega)\) to the unique solution \(v\) of scalar Chern-Simons problem in \(B_r(x_i)\) with datum \(\alpha\mu\).

Since \(v \le u\) and since points have zero \(W^{1, 2}\) capacity in \(\R^2\), by the Inverse maximum principle~\citelist{\cite{DP}*{Theorem~3} \cite{Ponce}*{Chapter~5}} we have for every \(i \in \{1, \dots, l\}\),
\[
-\Delta v(\{x_i\}) \le -\Delta u(\{x_i\}).
\]
Computing in particular both measure in the set \(\{x_i\}\), we get
\[
\alpha \mu(\{x_i\}) \le (\mu - \tau)(\{x_i\}) = \mu^\#(\{x_i\}),
\]
where
\[
\mu^\# = \mu -\tau.
\]
Taking the supremum over \(\alpha\), we deduce that
\[
2\pi \le \mu^\#(\{x_i\}).
\]
On the other hand, by Vázquez's nonexistence result \citelist{\cite{V}*{Section~5} \cite{BLOP}*{Section~5}}, we also have $\mu^\#(\{x_i\}) \le 2\pi$.
We conclude that
\[
\mu = \mu^\# \quad \text{in } \Omega \setminus\{x_1, \dots, x_l\}
\]
and for every \(i \in \{1, \dots, l\}\),
\[
\mu^\#(\{x_i\}) = 2\pi.
\]
In particular, the measure \(\mu^\#\) does not depend on the subsequence \((u_{n_k})_{k \in \N}\).
Since the solution of the Chern-Simons problem is unique for nonnegative datum, we deduce that the entire sequence \((u_n)_{n \in \N}\) converges to $u$ in $L^1(\Omega)$. 
The proof of the theorem is complete.
\qed

\section{Proof of Theorem~\ref{teoprinsist}}
\label{Thm34}

We first show that 
\[
\mu^\# \le \mu.
\]

Recall \(u_n \in W_0^{1, 1}(\Omega)\) and that for every \(\varphi \in C_c^\infty({\Omega})\),
\[
\label{eq:Equationun}
- \int_\Omega u_n \Delta \varphi + \int_\Omega \e^{v_n} (\e^{u_n} - 1)\varphi 
= \int_\Omega \varphi \dif\mu_n.
\]
The nonlinear term in the equation verified by \(u_n\) satisfies the sign condition: for every \(t \in \R\),
\[
\e^{v_n} (\e^t - 1) \sign{t} \ge 0.
\]
From the comparison estimate \cite{Ponce}*{Corollary~7.9}, \(\mu_n \ge 0\) implies that \(u_n \ge 0\).
Since the sequence \((u_n)_{n \in \N}\) and \((v_n)_{n \in \N}\) converge to \(u\) and \(v\) in \(L^1(\Omega)\), 
if the test function satisfies \(\varphi \ge 0\), then by Fatou's lemma,
\[
\int_\Omega \e^{v} (\e^u - 1) \varphi 
\le \liminf_{n \to \infty}{\int_\Omega \e^{v_n} (\e^{u_n} - 1) \varphi}. 
\]
As we let \(n\) tend to infinity in \eqref{eq:Equationun}, we get
\[
\int_\Omega \varphi \dif\mu^\#
= - \int_\Omega u \Delta \varphi + \int_\Omega \e^{v} (\e^{u} - 1) \varphi
\le \int_\Omega \varphi \dif\mu.
\]
Since this property holds for every \(\varphi \in C_c^\infty({\Omega})\) such that  \(\varphi \ge 0\), we deduce that \(\mu^\# \le \mu\).

\medskip

We now show that
\[
\mu^\# \ge 0.
\]

This property is proved in \cite{MAP}*{Theorem~1.3} in the case of semilinear equations with nonlinearities without dependence on the domain variable.
We explain below the main steps of the argument in our case.

Given an increasing sequence of nonnegative integers \((n_k)_{k \in \N}\), write
\[
\e^{v_{n_k}} (\e^{u_{n_k}} - 1)
= \e^{v_{n_k}} (\e^{u_{n_k}} - 1) \chi_{A_k} + \e^{v_{n_k}} (\e^{u_{n_k}} - 1) \chi_{\Omega \setminus A_k},
\]
where
\[
A_k = \{u_k \le k\} \cap \{v_k \le k\}.
\]
Using Cantor's diagonal argument, the sequence \((n_k)_{k \in \N}\) may be chosen such that \cite{MAP}*{Lemma~3.2}
\[
\e^{v_{n_k}} (\e^{u_{n_k}} - 1) \chi_{A_k} \to \e^{v} (\e^{u} - 1)
\quad \text{in \(L^1(\Omega)\).}
\]
By the capacitary estimate satisfied by the functions \(u_k\) and \(v_k\) \citelist{\cite{MAP}*{Lemma~3.2} \cite{Ponce}*{Lemma~9.4}},
\[
\capt_{W^{1, 2}}{(\Omega \setminus A_k)}
\le \capt_{W^{1, 2}}{(\{u_k > k\})} + \capt_{W^{1, 2}}{(\{v_k > k\})}
\le \frac{C}{k},
\]
for some constant \(C > 0\) independent of \(k\).
In particular, the \(W^{1, 2}\) capacity of the set \(\Omega \setminus A_k\) converges to zero as \(k\) tends to infinity.
Therefore, \(\big(\e^{v_{n_k}} (\e^{u_{n_k}} - 1) \chi_{\Omega \setminus A_k}\big)_{n \in \N}\) is a concentrating sequence with respect to the \(W^{1, 2}\) capacity in the sense of the Biting lemma \citelist{\cite{MAP}*{Section~2} \cite{BroCha:80}}.

We may now proceed as in the proof of \cite{MAP}*{Theorem~5.1} to conclude that \(\mu^\# \ge 0\).
The main ingredient in this step is a counterpart of the Inverse maximum principle concerning the concentrated limit of  sequences with respect to the \(W^{1, 2}\)~capacity \cite{MAP}*{Theorem~4.2}.

\medskip

We have proved that \(0 \le \mu^\# \le \mu\).
Reverting the roles of \(u\) and \(v\), we obtain \(0 \le \nu^\# \le \nu\).

\medskip
It remains to establish Assertion~\((iii)\).
For this purpose, let $\tau_1$ and $\tau_2$ be finite measures such that
\begin{equation}
\label{eqsm}
 \begin{aligned}
  \begin{aligned}
   \e^{v_{n}}(\e^{u_{n}}-1)&\stackrel{\ast}{\rightharpoonup}\e^v(\e^u-1)+\tau_1\\
   \e^{u_{n}}(\e^{v_{n}}-1)&\stackrel{\ast}{\rightharpoonup}\e^u(\e^v-1)+\tau_2
  \end{aligned}
  \quad \mbox{in } \mathcal{M}(\Omega).
 \end{aligned}
\end{equation}
Thus, $(u, v)$ solves the Chern-Simons problem
 \begin{equation*}
  \label{sistconv2}
   \left\{
    \begin{alignedat}{2}
      -\Delta u+\e^v(\e^u-1)&=\mu-\tau_1
       && \quad \mbox{in }\Omega,\\
      -\Delta v+\e^u(\e^v-1)&=\nu-\tau_2 
       && \quad \mbox{in }\Omega,\\
       u=v& =0
       && \quad \mbox{on }\partial \Omega.
    \end{alignedat}                      
\right.
  \end{equation*}
Proceeding as in the previous proof, we can apply the Brezis-Merle inequality and comparison estimates to show that $\tau_1$ and $\tau_2$ are supported in the finite set 
\[
B = \big\{x\in \Omega : \mu(\{x\})+\nu(\{x\})\geq 4\pi \big\}.
\]

Assuming that \(B\) is nonempty, we may write
\[
B = \{y_1, \dots, y_l\}
\]
where the points \(y_i \in \Omega\) are distinct.
Given \(i \in \{1, \dots, l\}\), let \(r > 0\) be such that \(B_r(y_i) \cap B = \{y_i\}\).
Adding the equations satisfied by \(u_n\) and \(v_n\), we have
\[
-\Delta (u_n+v_n)+2(\e^{u_n+v_n}-1) \geq \mu_n+\nu_n-2 \quad \mbox{in }\Omega.
\]
Note that for every 
\[
0 \le \alpha < \frac{4\pi}{\mu(\{y_i\}) + \nu(\{y_i\})},
\]
there exists \(w_n\) satisfying the equation
 \begin{equation*}
  \left\{
   \begin{alignedat}{2}
    -\Delta w_n+ 2(\e^{w_n}-1)&= \alpha(\mu_n+\nu_n)-2    && \quad\mbox{in } B_r(y_i),\\
      w_n&=0 && \quad \mbox{on } \partial B_r(y_i).
   \end{alignedat}
  \right.
 \end{equation*}
The existence of $w_{n}$ follows from \citelist{\cite{V}*{Theorem~2} \cite{Ponce}*{Chapter~8}}; alternatively, one may apply the method of sub and supersolution~\citelist{\cite{MP}*{Corollary~5.4} \cite{Ponce}*{Chapter~6}} with subsolution \(0\) and supersolution \(u_{n} + v_n\). 
In particular,
\[
0 \le w_n \le u_{n} + v_n \quad \text{in } \Omega.
\]

By a variant of Theorem~\ref{teocc} with nonlinearity \(\e^t(\e^t - 1)\) replaced by \(\e^t - 1\), the sequence \((w_n)_{n \in \N}\) converges in $L^1(B_r(y_i))$ to the solution of
\begin{equation*}
  \left\{
   \begin{alignedat}{2}
    -\Delta w+2(\e^{w}-1)&= \alpha(\mu+\nu)-2 &&    \quad \mbox{in } B_r(y_i),\\
                       w&=0 &&\quad \mbox{on } \partial B_r(y_i).
   \end{alignedat}
  \right.
 \end{equation*}
In particular, $w\leq u+v$ in \(B_r(y_i)\). 
By the Inverse maximum principle~\citelist{\cite{DP}*{Theorem~3} \cite{Ponce}*{Chapter~5}}, we deduce that
 \begin{equation*}
\alpha (\mu(\{y_i\}) + \nu(\{y_i\}))
\le \mu^\#(\{y_i\}) + \nu^\#(\{y_i\}). 
 \end{equation*}
Taking the supremum over \(\alpha\),
we conclude that
 \begin{equation*}
4\pi \le \mu^\#(\{y_i\}) + \nu^\#(\{y_i\}). 
 \end{equation*}
Since the reverse inequality holds, equality follows for every \(i \in \{1, \dots, l\}\).
The proof of the theorem is complete.
\qed

\section{Proof of Theorem~\ref{thmChernSimonsSystemAtomicPart}}
\label{Thm34-bis}

If \(\mu(\{a\}) + \nu(\{a\}) \le 4\pi\), then by Theorem~\ref{teoprinsist},
\[
\mu^\#(\{a\}) + \nu^\#(\{a\}) = \mu(\{a\}) + \nu(\{a\}).
\]
Since \(\mu^\# \le \mu\) and \(\nu^\# \le \nu\), we deduce that
\[
\mu^\#(\{a\}) = \mu(\{a\}) 
\quad\text{and}\quad
\nu^\#(\{a\}) = \nu(\{a\}).
\]

\medskip
We now assume that \(\mu(\{a\}) + \nu(\{a\}) > 4\pi\). 
In this case, by Theorem~\ref{teoprinsist},
\begin{equation}
\label{eqAssumption}
\mu^\#(\{a\}) + \nu^\#(\{a\}) = 4\pi.
\end{equation}

Recall that \(0 \le \mu^\# \le \mu\).
Thus, if \(\mu(\{a\}) = 0\), then \(\mu^\#(\{a\}) = 0\), whence \(\nu^\#(\{a\}) = 4\pi\) by the above identity.
Similarly, if \(\nu(\{a\}) = 0\), then \(\nu^\#(\{a\}) = 0\) and \(\mu^\#(\{a\}) = 4\pi\).
This concludes the proof of Assertion~\((i)\).

\medskip
In order to complete the proof of Assertion~\((ii)\), we assume that in addition to \eqref{eqAssumption}, we have
\begin{equation}
\label{eqAssumption-bis}
\mu(\{a\})\leq 4\pi
\quad\text{and}\quad
\nu(\{a\})\leq 4\pi.
\end{equation}
Using the notation of the proof of Theorem~\ref{teoprinsist}, we show that
\[
\tau_1(\{a\}) = \tau_2(\{a\}).
\]

Since $v_{n}\geq 0$ in \(\Omega\),
\begin{equation*}
\e^{v_{n}}(\e^{u_{n}}-1)\geq \e^{u_{n}}-1 \quad \mbox{in }\Omega.
 \end{equation*}
In particular, \(u_{n}\) is a subsolution of the Dirichlet problem
\begin{equation}
\label{eqDirichletScalaire}
  \left\{
   \begin{alignedat}{2}
-\Delta w + \e^{w}-1&= \lambda  &&\quad \mbox{in } \Omega,\\
w & =0        &&\quad    \mbox{on } \partial \Omega.
   \end{alignedat}
  \right.
 \end{equation}
with datum \(\lambda = \mu_n\).
Since \(\nu_n \ge 0\) and for every \(x \in \Omega\), \(\mu_n(\{x\}) + \nu_n(\{x\}) \le 4\pi\), we have for every \(x \in \Omega\), \(\mu_n(\{x\}) \le 4\pi\).
By V\'azquez existence result \citelist{\cite{V}*{Theorem~2} \cite{Ponce}*{Chapter~8}}, there exists \(\overline{u}_{n}\) satisfying the Dirichlet problem above with datum \(\lambda = \mu_n\).
By a comparison principle between the subsolution and the solution of the Dirichlet problem \cite{Ponce}*{Chapter~5}, $u_{n}\leq \overline{u}_{n}$ in $\Omega$. 

It follows from a variant of Theorem~\ref{teocc} with nonlinearity \(\e^t(\e^t - 1)\) replaced by \(\e^t - 1\) that the sequence \((\overline u_n)_{n \in \N}\) converges in \(L^1(\Omega)\) to the function \(\overline{u}\) satisfying the Dirichlet problem \eqref{eqDirichletScalaire} with datum \(\lambda = \tilde\mu\), where \(\tilde{\mu}\) is the largest measure less than or equal to \(\mu\) such that for every \(x \in \Omega\),
\[
\tilde\mu(\{x\}) \le 4\pi.
\]
In particular, since \(\mu(\{a\}) \le 4\pi\), we have
\[
\tilde\mu(\{a\}) = \mu(\{a\}).
\]
We also observe that the measure \(\mu - \tilde \mu\) is supported in a finite set, thus there exists \(r_1 > 0\) such that
\[
\tilde{\mu} = \mu \quad \text{in \(B_{r_1}(a)\).}
\]
Hence,
\begin{equation}\label{eqWeakConvergence1}
\e^{\overline u_{n}}\stackrel{\ast}{\rightharpoonup}\e^{\overline u} \quad \mbox{in } \mathcal{M}(B_{r_1}(a)).
\end{equation}

Similarly, if \(\overline{v}_n\) denotes the solution of the Dirichlet problem ~\eqref{eqDirichletScalaire} with datum \(\lambda = \nu_n\), then \(v_n \le \overline{v}_n\)  in \(\Omega\) and the sequence \((\overline{v}_n)_{n \in \N}\) converges to the solution of the Dirichlet problem \eqref{eqDirichletScalaire} with datum \(\lambda = \tilde\nu\) where the measure \(\tilde{\nu}\) satisfies \(\tilde{\nu}(\{a\}) = {\nu}(\{a\})\) and \(\nu - \tilde{\nu}\) is supported in a finite set.
In particular, there exists \(r_2 > 0\) such that
\begin{equation}\label{eqWeakConvergence2}
\e^{\overline v_{n}}\stackrel{\ast}{\rightharpoonup}\e^{\overline v} \quad \mbox{in } \mathcal{M}(B_{r_2}(a)).
\end{equation}

On the other hand, writing
\[
\e^{u_{n}}-\e^{v_{n}} = - \e^{u_{n}}(\e^{v_{n}}-1) + \e^{v_{n}}(\e^{u_{n}}-1),
\]
it follows from \eqref{eqsm} that
 \[
  \e^{u_{n}}-\e^{v_{n}}\stackrel{\ast}{\rightharpoonup}\e^u-\e^v+\tau_1-\tau_2 \quad \text{in } \cM(\Omega).
 \]
We observe that for every \(n \in \N\),
\[
-\e^{\overline v_n} \le \e^{u_{n}}-\e^{v_{n}} \le \e^{\overline u_n}
\quad \text{in \(\Omega\).}
\] 
As we let \(n\) tend to infinity, we deduce from \eqref{eqWeakConvergence1} and \eqref{eqWeakConvergence2} that for every \(0 < r \le \min{\{r_1, r_2\}}\),
\[
-\e^{\overline v} \le \e^{u}-\e^{v} + \tau_1-\tau_2 \le \e^{\overline u}
\quad \text{in \(B_r(a)\)}.
\]
Since the measure \(\tau_1 - \tau_2\) is supported in a finite set --- in particular is singular with respect to the Lebesgue measure --- we conclude that \(\tau_1 = \tau_2\).

\medskip
Let \(\tau = \tau_1 = \tau_2\).
By Theorem~\ref{teoprinsist}, we have
\[
\mu(\{a\}) + \nu(\{a\}) - 2 \tau(\{a\}) = \mu^\#(\{a\}) + \nu^\#(\{a\}) = 4\pi.
\]
Thus,
\[
\tau(\{a\}) =  \frac{\mu(\{y_i\}) + \nu(\{y_i\}) - 4\pi}{2},
\]
from which the conclusion follows.
\qed

\section{Concluding remarks}
\label{fr}

\subsection{Connection to reduced limits}
We have restricted ourselves to solutions of the Dirichlet problem, but we could also ask what happens to nonnegative solutions of the scalar Chern-Simons equation
\[
-\Delta u + \e^u(\e^u - 1) = \mu \quad \mbox{in } \Omega,
\]
without taking into account the Dirichlet boundary condition.

Solutions of the equation depend on the boundary data, but the approximation scheme does not.
More precisely, for each \(n \in \N\) take a solution \(u_n\) of the equation with datum \(\mu_n\) without prescribing any boundary condition. 
If the sequence \((u_n)_{n \in \N}\) converges to a function \(u\) in \(L^1(\Omega)\) and if the sequence \((\mu_n)_{n \in \N}\) converges weakly to some measure \(\mu\), one shows that \(u\) satisfies the scalar Chern-Simons equation with some datum \(\mu^\#\), possibly different from \(\mu\).
If we now take another sequence of solutions \((v_n)_{n \in \N}\) with the same data \((\mu_n)_{n \in \N}\) 
converging to another function \(v\) in \(L^1(\Omega)\), then \(v\) satisfies the scalar Chern-Simons equation with the same datum \(\mu^\#\).
This remarkable property has been recently discovered by Marcus and Ponce~\cite{MAP}, where they introduce the concept of reduced limit  \(\mu^\#\).

Combining Theorem~\ref{teocc} with \cite{MAP}*{Theorem~1.2}, we deduce the following result.

\begin{corolario}
Let $(\mu_n)_{n \in \N}$ be a nonnegative sequence of measures such that for every \(n \in \N\) and for every \(x \in \Omega\), 
\[
\mu_n(\{x\}) \le 2\pi
\]
and let \(u_n\) satisfy the scalar Chern-Simons problem
\[
-\Delta u_n + \e^{u_n}(\e^{u_n} - 1)  = \mu_n \quad \mbox{in } \Omega.
\]
If the sequence \((\mu_n)_{n \in \N}\) converges to a measure \(\mu\) in the sense of measures in \(\Omega\) and if the sequence \((u_n)_{n \in \N}\) converges to \(u\) in \(L^1(\Omega)\),
then \(u\) is the solution of the scalar Chern-Simons equation with datum \(\mu^\#\) defined in Theorem~\ref{teocc}. 
\end{corolario}

\subsection{Signed measures}
The sign of the measure $\mu$ affects substantially the conclusion. 
If \((\mu_n)_{n \in \N}\) is any sequence of nonpositive measures converging weakly in measure to some measure \(\mu\), then \(\mu\) is nonpositive and the sequence of solutions of the Dirichlet problem for the scalar Chern-Simons equation converge to the solution with datum \(\mu\).
This case is easier since the solutions \(u_n\) are nonpositive, whence the nonlinear term of exponential type is harmless.

The situation is more delicate when the sequence \((\mu_n)_{n \in \N}\) is not assumed to have a fixed sign.
In this case, one can show that given a signed measure $\mu$, positive numbers $c_1,\ldots, c_m$, and points $x_1,\dots, x_m\in \Omega$, there exists a sequence $(f_n)_{n \in \N}$ in $C_c^\infty(\Omega)$
such that
\begin{enumerate}[\((a)\)]
\item $(f_n)_{n \in \N}$ is bounded in \(L^1(\Omega)\),
\item $(f_n)_{n \in \N}$ converges to \(\mu\) in the sense of measures,
\item the solutions $u_n$ of the scalar Chern-Simons problem with datum \(f_n\) converge in \(L^1(\Omega)\) to a solution of the scalar Chern-Simons problem with datum \(\mu-\sum\limits_{i=1}^m c_i\delta_{x_i}\).
\end{enumerate}

In particular, it is not possible in this case to have an explicit formula of the measure \(\mu^\#\) only in terms of \(\mu\). We refer to \cite{Ad}*{Teorema 4.8} for the proof. 

\subsection{Nonuniqueness of the reduced limit of the Chern-Simons system}
\label{NRDCSS}

The reduced limit \((\mu^\#, \nu^\#)\) of the Chern-Simons system cannot be computed only in terms of the weak\(^*\) limit \((\mu, \nu)\) when both conditions
\begin{enumerate}[\((a)\)]
\item \(\mu(\{x\}) = 0\) or \(\nu(\{x\}) = 0\),
\item \(\mu(\{x\})\leq 4\pi\) and \(\nu(\{x\})\leq 4\pi\),
\end{enumerate}
fail for some \(x \in \Omega\).

\begin{proposicao}
Let \(a \in \Omega\). 
For every \(\alpha > 4\pi\) and for every \(\beta > 0\), there exist sequences \((\mu_n^i)_{n \in \N}\) and \((\nu_n^i)_{n \in \N}\) of nonnegative functions in \(L^1(\Omega)\) with \(i \in \{1, 2\}\) such that
\begin{enumerate}[\((i)\)]
\item \((\mu_n^i)_{n \in \N}\) and \((\nu_n^i)_{n \in \N}\) converge weakly to \(\alpha\delta_a\) and \(\beta\delta_a\) in the sense of measures in \(\Omega\),
\item there exists a sequence of solutions \(\big( (u_n^i, v_n^i ) \big)_{n \in \N}\) of the Chern-Simons system with datum \((\mu_n^i, \nu_n^i)\) converging in \(L^1(\Omega) \times L^1(\Omega)\) to a solution with datum \((\alpha^{i, \#} \delta_a , \beta^{i, \#} \delta_a)\),
\item \(\alpha^{1, \#} \ne \alpha^{2, \#}\) and \(\beta^{1, \#} \ne \beta^{2, \#}\).
\end{enumerate}
\end{proposicao}

We shall not prove this proposition.
Instead, we restrict ourselves to the case where \((\alpha\delta_a, \beta\delta_a)\) is given by
\[
(5\pi \delta_a, 2\pi \delta_a)
\]
for some \(a \in \Omega\) in order to emphasize the main idea of the proof.

For this purpose, let \((f_n)_{n \in \N}\) be a sequence of nonnegative functions in \(L^1(\Omega)\) such that
\[
f_n \overset{*}{\rightharpoonup} \delta_a
\quad \text{in \(\cM(\Omega)\).}
\]

We construct the first sequence \(\big((\mu_n^1, \nu_n^1)\big)_{n \in \N}\) of the form
\[
\big((5\pi f_{m_n}, 2\pi f_n)\big)_{n \in \N}
\]
where \((m_n)_{n \in \N}\) is a sequence of positive integers to be chosen below.
For fixed \(n \in \N\), let \(\big((u_{m, n}^1, v_{m, n}^1)\big)_{m \in \N}\) be a solution of the Chern-Simons system with datum \((5\pi f_{m}, 2\pi f_n)\).
Then,  as \(m\) tends to infinity,
\[
(5\pi f_{m}, 2\pi f_n) \overset{*}{\rightharpoonup} (5\pi \delta_a, 2\pi f_n)  \quad \text{in \(\cM(\Omega) \times \cM(\Omega)\).}
\]

By a standard property of elliptic equations with absorption term \cite{Ponce}*{Chapter~7} and by Stampacchia's linear regularity theory~\cite{Ponce}*{Chapter~3}, the sequence 
\(\big((u_{m, n}^1, v_{m, n}^1)\big)_{m \in \N}\) is compact in \(L^1(\Omega) \times L^1(\Omega)\).
It is then possible to extract a subsequence with respect to the index \(m\) if necessary such that for every \(n \in \N\), 
\[
(u_{m, n}^1, v_{m, n}^1) \to (u_n^1, v_n^1)
\quad
\text{in  \(L^1(\Omega) \times L^1(\Omega)\)}
\]
as \(m\) tends to infinity.
It follows from Theorem~\ref{thmChernSimonsSystemAtomicPart} that \((u_n^1, v_n^1)\) satisfies the Chern-Simons system with datum \((4\pi \delta_a, 2\pi f_n)\).

Note that 
\[
(4\pi \delta_a, 2\pi f_n) \overset{*}{\rightharpoonup} (4\pi \delta_a, 2\pi \delta_a)  \quad \text{in \(\cM(\Omega) \times \cM(\Omega)\).}
\]
By compactness of the sequence \(\big((u_n^1, v_n^1)\big)_{n \in \N}\) in \(L^1(\Omega) \times L^1(\Omega)\), we may extract a subsequence converging to \((u^1, v^1)\).
By Theorem~\ref{thmChernSimonsSystemAtomicPart}, \((u^1, v^1)\) satisfies the Chern-Simons system with datum \((3\pi \delta_a, \pi \delta_a)\).

For every \(n \in \N\), take \(m_n \in \N\) such that
\[
\|u_{m_n, n}^1 - u_n^1\|_{L^1(\Omega)} + \|v_{m_n, n}^1 - v_n^1\|_{L^1(\Omega)}
\le \frac{1}{n + 1}.
\]
Then, the sequence \(\big((u_{m_n, n}^1, v_{m_n, n}^1)\big)_{m \in \N}\) converges to \((u^1, v^1)\) in \(L^1(\Omega) \times L^1(\Omega)\).
We have found a sequence of solutions of the Chern-Simons system with datum \(\big((5\pi f_{m_n}, 2\pi f_n)\big)_{n \in \N}\) converging to the solution with datum 
\[
\boxed{(3\pi \delta_a, \pi \delta_a)}.
\]

\medskip
We construct the second sequence \((\mu_n^2, \nu_n^2)_{n \in \N}\) of the form
\[
\big((4\pi f_{m_n} + \pi f_n, 2\pi f_{m_n})\big)_{n \in \N},
\]
where \((m_n)_{n \in \N}\) is a sequence of positive integers, possibly different from the previous one.
For fixed \(n \in \N\), let \(\big((u_{m, n}^2, v_{m, n}^2)\big)_{m \in \N}\) be a solution of the Chern-Simons system with datum \((4\pi f_{m} + \pi f_n, 2\pi f_m)\).
Then,  as m tends to infinity,
\[
(4\pi f_{m} + \pi f_n, 2\pi f_m) \overset{*}{\rightharpoonup} (4\pi \delta_a + \pi f_n, 2\pi \delta_a)  \quad \text{in \(\cM(\Omega) \times \cM(\Omega)\).}
\]
The sequence \(\big((u_{m, n}^2, v_{m, n}^2)\big)_{m \in \N}\) is compact in \(L^1(\Omega) \times L^1(\Omega)\).
It is then possible to extract a subsequence of with respect to the index \(m\) if necessary such that for every \(n \in \N\), 
\[
(u_{m, n}^2, v_{m, n}^2) \to (u_n^2, v_n^2)
\quad
\text{in  \(L^1(\Omega) \times L^1(\Omega)\).}
\]
It follows from Theorem~\ref{thmChernSimonsSystemAtomicPart} that \((u_n^2, v_n^2)\) satisfies the Chern-Simons system with datum \((3\pi \delta_a + \pi f_n, \pi \delta_a)\).

Note that 
\[
(3\pi \delta_a + \pi f_n, \pi \delta_a) \overset{*}{\rightharpoonup} (4\pi \delta_a, \pi \delta_a)  \quad \text{in \(\cM(\Omega) \times \cM(\Omega)\).}
\]
By compactness of the sequence \(\big((u_n^2, v_n^2)\big)_{n \in \N}\) in \(L^1(\Omega) \times L^1(\Omega)\), we may extract a subsequence converging to \((u^2, v^2)\).
By Theorem~\ref{thmChernSimonsSystemAtomicPart}, \((u^2, v^2)\) satisfies the Chern-Simons system with datum \((\frac{7\pi}{2} \delta_a, \frac{\pi}{2} \delta_a)\).

Proceeding as before, for every \(n \in \N\) we may choose \(m_n \in \N\) such that
\(\big((u_{m_n, n}^2, v_{m_n, n}^2)\big)_{m \in \N}\) converges to \((u^2, v^2)\) in \(L^1(\Omega) \times L^1(\Omega)\).
Hence, there exists a sequence of solutions of the Chern-Simons system with datum \(\big((4\pi f_{m_n} + \pi f_n, 2\pi f_{m_n})\big)_{n \in \N}\) converging to the solution with datum
\[
\boxed{\big(\tfrac{7\pi}{2} \delta_a, \tfrac{\pi}{2} \delta_a\big)}.
\]

\section*{Acknowledgments}
The first author (ACP) was supported by the Fonds de la Recherche scientifique---FNRS. 
The second author (AEP) thanks the Brazilian research agencies CAPES and CNPq under grants 446809-0 and 142584/2008-8, respectively, for the financial support and the Université catholique de Louvain for its hospitality during the development of part of this work. 
%
%
%
%



\begin{bibdiv}
\begin{biblist}

\bib{BLOP}{article}{
      author={Bartolucci, Daniele},
      author={Leoni, Fabiana},
      author={Orsina, Luigi},
      author={Ponce, Augusto~C.},
       title={Semilinear equations with exponential nonlinearity and measure
  data},
        date={2005},
     journal={Ann. Inst. H. Poincar\'e Anal. Non Lin\'eaire},
      volume={22},
       pages={799\ndash 815},
}

\bib{BenBre:04}{article}{
      author={B{\'e}nilan, {\mbox{Ph}}ilippe},
      author={Brezis, Ha\"{\i}m},
       title={Nonlinear problems related to the {T}homas-{F}ermi equation},
        date={2004},
     journal={J. Evol. Equ.},
      volume={3},
       pages={673\ndash 770},
        note={Dedicated to Ph.~B\'enilan},
}

\bib{Bre:80}{article}{
      author={Brezis, Ha{\"{\i}}m},
       title={Some variational problems of the {T}homas-{F}ermi type},
conference={
      title={Variational inequalities and complementarity problems ({P}roc.
        {I}nternat. {S}chool, {E}rice, 1978)},
   },
   book={
      publisher={Wiley},
      place={Chichester},
      date={1980},
   },
       pages={53\ndash 73},
}

\bib{Bre:82}{article}{
      author={Brezis, Ha{\"{\i}}m},
       title={Probl\`emes elliptiques et paraboliques non lin\'eaires avec
  donn\'ees mesures},
conference={
      title={Goulaouic-{M}eyer-{S}chwartz {S}eminar, 1981/1982},         
      },
   book={
      publisher={\'Ecole Polytech.},
      place={Palaiseau},
      date={1982},
   },
  note={Exp. {N}o. {XX}, 13},
}

\bib{BS}{article}{
      author={Brezis, Ha{\"{\i}}m},
      author={Strauss, Walter~A.},
       title={Semi-linear second-order elliptic equations in {$L^{1}$}},
        date={1973},
     journal={J. Math. Soc. Japan},
      volume={25},
       pages={565\ndash 590},
}




\bib{BPM}{article}{
      author={Brezis, Haïm},
      author={Marcus, M.},
      author={Ponce, Augusto~C.},
      title={Nonlinear elliptic equations with measures revisited},
conference={
      title={Mathematical aspects of nonlinear dispersive equations},
   },
   book={
      editor={{B}ourgain, {J}.},
      editor={{K}enig, {C}.},
      editor={{K}lainerman, {S}.},
      series={Ann. of Math. Stud.},
      volume={163},
      publisher={Princeton Univ. Press},
      place={Princeton, NJ},
   },
     PAGES = {55--110},
}

\bib{BM}{article}{
      author={Brezis, Haïm},
      author={Merle, Frank},
       title={Uniform estimates and blow-up behavior for solutions of
  ${-\Delta} u={V}(x)e^u$ in two dimensions},
        date={1991},
     journal={Comm. Partial Differential Equations},
      number={16},
       pages={1223\ndash 1253},
}

\bib{BroCha:80}{article}{
      author={Brooks, J.~K.},
      author={Chacon, R.~V.},
       title={Continuity and compactness of measures},
        date={1980},
     journal={Adv. in Math.},
      volume={37},
       pages={16\ndash 26},
}

\bib{DP}{article}{
      author={Dupaigne, Louis},
      author={Ponce, Augusto~C.},
       title={Singularities of positive supersolutions in elliptic {PDEs}},
        date={2004},
     journal={Selecta Math. (N.S.)},
      volume={10},
       pages={341\ndash 358},
}

\bib{Evans_Gariepy}{book}{
   author={Evans, Lawrence C.},
   author={Gariepy, Ronald F.},
   title={Measure theory and fine properties of functions},
   series={Studies in Advanced Mathematics},
   publisher={CRC Press},
   place={Boca Raton, FL},
   date={1992},
}

\bib{LPY}{article}{
      author={Lin, Chang-Shou},
      author={Ponce, Augusto~C.},
      author={Yang, Yisong},
       title={A system of elliptic equations arising in {C}hern-{S}imons field
  theory},
        date={2007},
     journal={J. Funct. Anal.},
      volume={247},
       pages={289\ndash 350},
}

\bib{MAP}{article}{
      author={Marcus, M.},
      author={Ponce, Augusto~C.},
       title={Reduced limits for nonlinear equations with measures},
        date={2010},
     journal={J. Funct. Anal.},
      number={258},
       pages={2316\ndash 2372},
}

\bib{MP}{article}{
      author={Montenegro, Marcelo},
      author={Ponce, Augusto~C.},
       title={The sub-supersolution method for weak solution},
        date={2008},
     journal={Proc. Amer. Math. Soc.},
      volume={136},
       pages={2429\ndash 2438},
}

\bib{Ponce}{misc}{
      author={Ponce, Augusto~C.},
       title={Selected problems on elliptic equations involving measures},
        date={2012},
        note={Available at \texttt{http://arxiv.org/abs/1204.0668}},
}

\bib{Ad}{thesis}{
      author={Presoto, Adilson~E.},
       title={Soluções limites para problemas elípticos envolvendo medidas},
        type={Ph.D. Thesis},
       place={IMECC---UNICAMP, Campinas, SP, Brazil}
        date={2011},
}

\bib{Stampacchia}{article}{
      author={Stampacchia, Guido},
       title={Le probl\`eme de {D}irichlet pour les \'equations elliptiques du
  second ordre \`a coefficients discontinus},
        date={1965},
     journal={Ann. Inst. Fourier (Grenoble)},
      volume={15},
       pages={189\ndash 258},
}

\bib{V}{article}{
      author={V{\'a}zquez, Juan~L.},
       title={On a semilinear equation in {${\bf R}^{2}$} involving bounded
  measures},
        date={1983},
     journal={Proc. Roy. Soc. Edinburgh Sect. A},
      volume={95},
       pages={181\ndash 202},
}

\end{biblist}
\end{bibdiv}

\end{document}